# Why There Are No Gaps In The Support Of
# Non-Negative Integer-Valued Infinitely Divisible Laws?

**S Satheesh**

(*Bharat Sanchar Nigam Limited*)

NEELOLPALAM, S. N. Park Road
Trichur – 680 004, **India.**

ssatheesh@sancharent.in

**Abstract.** Remark.9 in Bose-Dasgupta-Rubin (2002) review states that when a non-negative integer-valued infinitely divisible law has an atom at unity then its support cannot have any gaps. Here one has two questions. (i) Why there are no gaps and (ii) Can there be gaps if the condition is not satisfied. Our investigation with these questions in mind centers on the implications of having and not having atoms at zero and unity. We give two examples/ constructions, which show that the remark needs modification and we modify it.

## 1. Introduction

The Bose-Dasgupta-Rubin (2002) review (abbreviated as BDR (2002)) extensively reviews the literature on infinitely divisible (ID) laws and processes with many illustrative examples and a list of Levy measures. This paper is inspired by their remark.9, which states that when a non-negative integer-valued ID law has a positive probability at 1 then its support cannot have any gaps. Such a statement leads to two questions; (i) Why there are no gaps in the support of such a r.v $X$ when P{$X$=1}>0 and (ii) Can there be gaps if the condition is not satisfied. Interestingly this investigation shows that $X$ can have gaps in its support even when P{$X$=1}>0 and the remark.9 of BDR (2002) needs modification for it to be true. As it turns out a related problem in this context is the implication of P{$X$=0}>0 for a non-negative integer-valued ID r.v $X$.

In fact, in the literature one comes across phrases like "an ID distribution {$p_i$ , i=0,1,2, ....} with $p_0$>0", see eg. Katti (1967), Steutel (1973, 1979). So a natural question is what is the effect of $p_0 > 0$ or $P\{X=0\} > 0$ in terms of the r.v $X$ and whether there are ID lattice r.vs $X$ for which $P\{X=0\}=0$. Feller (1968, p.290) while proving that ID laws on $I_0 = \{0,1,2, \ldots.\}$ must be compound Poisson (and conversely) states the condition $P\{X=0\} > 0$ as an analytical requirement in terms of the corresponding probability generating function (PGF). Steutel and van Harn (1979) assume this condition in the statement of the result itself. Bondesson (1981) discusses ID laws on $I_0$ and those obtained by truncating out the mass at zero. But whether they continue to be ID or not is not mentioned. Johnson, *et al.* (1992, p.352) proves that if $X$ is compound Poisson (*ie.* ID) and has finite mean then $P\{X=0\} > 0$ and a partial converse of this. BDR (2002) also does not discuss this point. We know that the geometric r.v $X$ on $I_1 = \{1,2,3, \ldots.\}$ is ID and $P\{X=0\}=0$.

Here first we give a probabilistic argument to prove that the condition $P\{X=0\} > 0$ is necessary for all ID laws on $I_0$ with integer-valued components. Among other implications it also explains why certain ID laws on $I_0$ cannot have any gaps. Our discussion brings out certain important implications of $P\{X=0\} > 0$ and $P\{X=1\} > 0$. Finally we give an alternative derivation of the result of Feller (1968, p.290) that ID laws on $I_0$ must be compound Poisson (and conversely). This paper is part of the review, Satheesh (2003).

## 2. Integer-Valued Infinitely Divisible Laws.

**Theorem 1** If a r.v $X$ on $I_0$ is ID with integer-valued components then $P\{X=0\} > 0$.

**Proof.** Under the assumptions on $X$, there exists i.i.d r.vs $\{X_{in}\}$ on $I_0$ such that

$$X \overset{d}{=} X_{1n} + \ldots. + X_{nn} \quad \text{for every } n \geq 1 \text{ integer,} \qquad (1)$$

where $\{X_{in}\}$, the components are integer-valued. Now assume the contrary that $P\{X=0\}=0$. Then $P\{X_{in} = 0\}=0$ for all $1 \leq i \leq n$ and every $n \geq 1$. Let $k > 0$ be the least integer such that $P\{X = k\} > 0$. For a given $n \geq 1$, let $r > 0$ is the least integer such that $P\{X_{in} = r\} > 0$. Then the minimum value assumed with positive probability by the RHS of (1) is $nr$ where as that of the LHS is $k$ and they are never equal when $n > k$. Hence when $n > k$ we cannot have the representation (1). The minimum values $nr$ and $k$ are equal for every $n \geq 1$ only if $k = 0 = r$. But $k = 0$ implies $r = 0$. Hence $P\{X = 0\} > 0$. Also $P\{X_{in} = 0\} > 0$ for all $i = 1, \ldots, n$.

The author has come to know that Kallenberg has priority to this result, with a different proof, see Grandall (1997, p.26)). We also have:

**Corollary 1** Let $Q(s)$ is the PGF of a r.v $X$ on $I_0$ that is ID with integer-valued components. Then P{$X$=0}>0 or equivalently $Q$(o)>0. Hence $Q$(s) never vanishes. This can be seen, as the integer-valued analogue of the fact that the characteristic function of ID laws does not have real zeroes.

**Corollary 2** For an ID r.v $X$ on $I_0$ , P{$X$=0}=0 implies that the components {$X_{in}$} of $X$ are no more supported by positive integers.

Arguing on lines similar to the proof of theorem.1 we have:

**Theorem 2** If a r.v $X$ on $I_0$ is ID with integer-valued components {$X_{in}$}, then P{$X>k$}>0 for any $k$>0 integer. In other words P{$X>k$}=0 for some positive integer $k$ is impossible.

**Corollary 3** The support of an ID law on $I_0$ cannot be a bounded subset of $I_0$. Recall the well-known property that the support of ID laws cannot be bounded. This is why the binomial law is not ID.

Hu, *et al*. (2004) in their lemma.1 have shown that if P{$X$=0}>0 then the support of a discrete ID r.v $X \geq 0$ coincides with its components {$X_{in}$} for every $n \geq 1$. Hu, *et al*. (2004) reports also that this result is in Sato (1999, Theorem.24.5, p.149) or Steutel and van Harn (2004, Corollary.8.3, p.111). Thus for any PGF $Q(s)$ that is ID, the PGF $P(s) = sQ(s)$ is not ID with integer-valued components. *Eg*. the geometric law on $I_1$ is not ID with integer-valued components though the geometric law on $I_0$ is ID with integer-valued components. But we know that the property of infinite divisibility is invariant under translation. Hence the geometric law on $I_1$ also is ID but its components are not integer-valued. Thus:

**Property 1** The property "infinite divisibility with integer-valued components" is not invariant under translation.

Now, can we distinguish the class of ID laws on $I_0$ whose components are integer-valued and those whose components are not? From theorem.1 we have:

**Corollary 4** If the support of an ID r.v $X$ on $I_0$ coincides with that of its components {$X_{in}$} for every $n \geq 1$ integer, then {$X_{in}$} are also integer-valued and hence P{$X$=0}>0.

**Theorem 3** Let $X$ be a non-negative integer-valued ID r.v. Then P$\{X=0\}>0$ iff the support of $X$ coincides with that of its components $\{X_{in}\}$ for every $n \geq 1$ integer.

**Proof.** The sufficiency part follows by corollary.4 and the necessity part by lemma.1 of Hu, *et al*. (2004).

Combining these ideas we have:

**Theorem 4** For a non-negative integer-valued ID r.v $X$ the following statements are equivalent.

(a) P$\{X=0\}>0$

(b) The support of $X$ coincides with that of its components $\{X_{in}\}$ for every $n \geq 1$ integer.

(c) $X$ has integer-valued components.

**Proof.** (a) implies (b) by lemma.1 of Hu, *et al*. (2004). (b) implies (c) is obvious since $X$ is integer-valued. Finally (c) implies (a) by theorem.1.

Thus the results in the literature on non-negative integer-valued ID r.vs $X$ with P$\{X=0\}>0$ is valid only for the class identified in theorem.3. By assuming P$\{X=0\}>0$, the support of a sum of r.vs and the components in the sum are the same and so a requirement for describing stability and self-decomposability is satisfied in Steutel and van Harn (1979).

**Remark 1** With reference to the Remark.9 of BDR (2002) one may ask; Why there are no gaps in the support of a non-negative integer-valued ID r.v $X$ with P$\{X=1\}>0$? If this ID r.v $X$ has integer-valued components we then have both P$\{X=0\}>0$ and P$\{X=1\}>0$. Hence for every $n \geq 1$ the components $X_{in}$ also satisfy P$\{X_{in}=0\}>0$ and P$\{X_{in}=1\}>0$. Consequently $X$ cannot have any gaps in its support as the summation (1) ensures that every non-negative integer carries a probability.

An interesting question now is; are there any ID r.vs on $I_0$ with integer-valued components and have gaps in its support? Of course, one expects that they cannot have a probability at $X=1$. The answer is in the following examples.

**Example 1** From the PGF of a negative binomial law consider the r.v $X$ with PGF

$$\{p/(1-qs^k)\}^t, \ p+q=1, \ k>1 \text{ integer and } t>0.$$

This is the PGF of the negative binomial-sum of r.vs degenerate at $k>1$ integer. Obviously $X$ is ID, P{$X=0$}>0 (hence the components are integer-valued), the support of $X$ has gaps as its atoms are $k-1$ integers apart, but P{$X=1$}=0 (as $k>1$). In general, consider a PGF $Q(s)$ that is ID, $Q(o)>0$ and consider the r.v $X$ having PGF $Q(s^k)$. Now P{$X=0$}>0, P{$X=1$}=0 and the atoms of $X$ are $k-1$ integers apart.

However, the following example shows that when we do not impose the condition that the ID r.v $X$ has integer-valued components then even P{$X=1$}>0 cannot not rule out the possibility of having gaps in its support.

**Example 2** As a next step from example.1, now consider the r.v $X$ with PGF

$$s\{p/(1-qs^k)\}^t, \ p+q=1, \ k>1 \text{ integer and } t>0.$$

Here $X$ is ID, P{$X=0$}=0 (hence the components are not integer-valued), P{$X=1$}>0 and $X$ has gaps in its support (its atoms are $k-1$ integers apart). In general, consider a PGF $Q(s)$ that is ID, $Q(o)>0$ and consider the r.v $X$ having PGF $sQ(s^k)$. Now P{$X=0$}=0, P{$X=1$}>0 and the atoms of $X$ are $k-1$ integers apart.

**Remark 2** Thus remark.9 of BDR (2002) holds good only for ID laws on $I_0$ with P{$X=0$}>0. The following modification of the remark.9 is now clear:

**Theorem 5** Let $X$ be an ID r.v on $I_0$ with P{$X=0$}>0 (or any of the equivalent statements in theorem.4). Then the support of $X$ cannot have any gaps if (and only if) P{$X=1$}>0.

Now we derive the result of Feller (1968, p.290) that ID laws on $I_0$ must be compound Poisson (and conversely), in two parts for their independent interest.

**Theorem 6** If an ID law on [0,∞) has an atom at the origin then it must be compound Poisson.

**Proof.** Let $\varphi(s) = e^{-\psi(s)}$ be the Laplace transform (LT) of the ID law (Feller, 1971, p.450). If a distribution with LT $\varphi$ has an atom at the origin then it is equivalent to $\varphi(\infty) > 0$. For an ID law this is reflected as $\psi(\infty) = \lambda < \infty$. Setting $\psi(s)/\lambda = F(s)$, we have $F(o) = 0$ and $F(\infty) = 1$. Further since $\psi(s)$ has completely monotone derivative (CMD), $F(s)$ is continuous and non-decreasing. Thus $F(s)$ is a d.f with CMD. Hence by Pillai and Sandhya (1990) $F(s)$ is the d.f of a mixture of exponential laws and hence $F(s) = 1-\alpha(s)$, where $\alpha(s)$ is a LT. Hence $\varphi(s) = exp\{-\lambda(1-\alpha(s))\}$ completing the proof.

**Corollary 5** A probability distribution on $[0,\infty)$ with LT $\varphi$ is compound Poisson iff $\theta\{-ln\varphi(x)\}$, $x>0$ is a d.f for some $\theta>0$.

**Theorem 7** (As integer-valued analogue of ID laws on $[0,\infty)$). ID laws on $I_0$ having integer-valued components (equivalently in terms of PGFs, $Q = \{Q_n\}^n$, for every $n\geq1$ integer) are Poisson compounds of non-negative integer-valued r.vs and conversely.

**Proof.** Readily follows from Theorems.1 & 6 and by noticing that the LT $\omega(s)$ in theorem.6 must now be that of a distribution on $I_0$. The converse is clear.

We have come to know that a different proof of this result is available in Ospina and Gerber (1987) and its multivariate extension by Sundt (2000). Steutel and van Harn (1979) has discussed discrete analogues of self-decomposable and stable laws while Satheesh and Sandhya (2004) discusses discrete analogues of self-decomposable, stable and semi-stable laws and their generalizations.

**Acknowledgement.** The author thank Prof. M Sreehari, M S University of Baroda, India, for the line of argument in theorem.1, Prof. A M Iksanov, National University, Kiev, Ukrain, for his suggestion recorded as theorem.3 and some preprints and the referee whose careful reading removed some ambiguities and for the suggestion to formulate property.1.